\documentclass {article}
\usepackage[utf8]{inputenc}
\usepackage[english]{babel}
\usepackage{amsmath}
\usepackage{amsthm}
\usepackage{amssymb}

\input xy
\xyoption{all}

\begin{document}

    \begin{center}
        \large Alexander Zakharov \par \vspace{0.3 cm}
        \bf \Large Intersecting free subgroups in free amalgamated products of two groups with normal finite amalgamated subgroup
    \end{center}
    \vspace {0.5 cm}

    We partly generalize the estimate for the rank of intersection of subgroups in free products of groups, proved earlier by S.V.Ivanov and W.Dicks, to the case of free amalgamated products of groups with normal finite amalgamated subgroup. We also prove that the obtained estimate is sharp and cannot be further improved when the amalgamated product contains an involution.

    \vspace {0.5 cm}
     \subsection* {1. Introduction.} \par
    Let first $G$ be a free group, $H_1$ and $H_2$ --- finitely generated subgroups in $G$.
    In 1954 Howson \cite {1} proved that in this case subgroup $H_1 \cap H_2$ is also finitely generated. Then in 1957 Hanna Neumann \cite {2} established the following estimate for the rank of intersection of subgroups in a free group  (Hanna Neumann inequality):
    \begin {equation}\label{eq1}
          \overline{r}(H_1 \cap H_2) \leqslant 2\overline{r}(H_1)\overline{r}(H_2),
    \end {equation}
    where $\overline{r}(H)= max \: (0, r(H)-1)$ is the reduced rank of subgroup $H$.

    \par
      S.V.Ivanov and W.Dicks \cite{3}, \cite{4}, \cite{8} generalized these results to the case when $G$ is a free product of groups. Below we consider only nontrivial free products. We call a subgroup of a free product of groups {\it factor-free},  if it intersects trivially with the conjugates to the factors of the free product. Factor-free subgroups are free, according to the Kurosh subgroup theorem \cite{7}. We need the following theorem proved in \cite{4}:

    \newtheorem*{theorem}{Theorem 1. (W.Dicks, S.V.Ivanov)}
    \begin {theorem}
        Suppose $ G=G_{1}*G_{2}$ is a free product of groups, and $H_{1}$, $H_{2}$ are factor-free subgroups of $G$  with finite ranks. Then the intersection $H_{1}\cap H_{2}$ also has finite rank and
    \begin {equation}\label{eq2}
     {\overline{r}}(H_1 \cap H_2) \leqslant 2\frac{q^{*}}{q^{*}-2}{\overline{r}}(H_1){\overline{r}}(H_2) \leqslant 6{\overline{r}}(H_1){\overline{r}}(H_2),
    \end {equation}
     where $q^{*}$ is the minimum of orders $> 2$ of subgroups of groups $G_{1}$, $G_{2}$, and $\frac{q^{*}}{q^{*}-2} = 1$ if $q^{*} = \infty$. In addition, the first estimate in (\ref{eq2}) is sharp and cannot be further improved whenever $G$ contains an involution (element of order 2) and $G \ncong \mathbb{Z}_{2}*\mathbb{Z}_{2}$.
    \end {theorem}

    It is easy to see that this theorem generalizes Hanna Neumann inequality (\ref{eq1}).
    \vspace {0.1 cm}
    Below we prove a further generalization of inequalities (\ref{eq1}) and (\ref {eq2}) to the case of free amalgamated product with normal finite amalgamated subgroup. We also consider the question when the obtained estimate is sharp.

    \vspace {0.1 cm}

     \newtheorem*{theorem1}{Theorem 2}
    \begin {theorem1}
     Suppose $G = G_{1}*_{T}G_{2}$ is an amalgamated free product, $T$ is normal (in $G$) finite, and $H_{1}$, $H_{2}$ are factor-free (and, therefore, free) subgroups of $G$ with finite ranks. Then the intersection $H_{1}\cap H_{2}$ also has finite rank, and
    \begin {equation}\label{eq3}
     {\overline{r}}(H_1 \cap H_2) \leqslant 2\frac{q_{f}^{*}}{q_{f}^{*}-2}|T| \cdot {\overline{r}}(H_1){\overline{r}}(H_2) \leqslant 6|T| \cdot {\overline{r}}(H_1){\overline{r}}(H_2),
    \end {equation}
     where $q_{f}^{*}$ is the minimum of orders $> 2$ of subgroups of groups $G_{1}/T$, $G_{2}/T$, and $\frac{q_{f}^{*}}{q_{f}^{*}-2} = 1$ if $q_{f}^{*} = \infty$,  $|T|$ is the order of group $T$. In addition, the first estimate in (\ref{eq3}) is sharp and cannot be further improved whenever $G_{1}/T$ or $G_{2}/T$ contains an involution and $G_{1}/T *G_{2}/T \ncong \mathbb{Z}_{2}*\mathbb{Z}_{2}$.

    \end {theorem1}

    \subsection* {2. Proof of the estimate.}
      Here we prove the estimate (\ref{eq3}).
    \par
    \vspace{0.2 cm}

      Since $T$ is a normal subgroup of $G_{1}*_{T}G_{2}$,
    we can consider a factorization $$\varphi : G_{1}*_{T}G_{2} \rightarrow G_{1}/T*G_{2}/T.$$
    Let $$\varphi(H_{1})=H_{1}',\quad  \varphi(H_{2})=H_{2}', \quad
     \varphi(H_{1}\cap H_{2})=L\subseteq H_{1}'\cap H_{2}'$$
     The last inclusion holds since $$L =  \varphi(H_{1}\cap H_{2}) \subseteq \varphi (H_{1}) \cap \varphi (H_{2})=H_{1}'\cap H_{2}'$$

     \newtheorem{lemmas}{Lemma}
    \begin {lemmas} \label {lemma 1}
        \begin {equation} \label{eq4}
            \overline{r}(H_{1}\cap H_{2}) \leq 2\frac{q_{f}^{*}}{q_{f}^{*}-2}|H_{1}'\cap H_{2}':L|\: \overline{r}(H_{1})\overline{r}(H_{2})
          \end {equation}
    \end {lemmas}
    \par
    (It will be shown below, in lemma 2, that the index $|H_{1}'\cap H_{2}':L|$ is finite.) \par
    $\square$
    Since $H_{1}$ and $H_{2}$ are factor-free, they intersect trivially with the amalgamated subgroup, therefore $$H_{1}\simeq H_{1}', \quad H_{2}\simeq H_{2}', \quad H_{1}\cap H_{2}\simeq L,$$ so \par $$\overline{r}(H_{1})=\overline{r}(H_{1}'), \quad \overline{r}(H_{2})=\overline{r}(H_{2}'), \quad \overline{r}(H_{1}\cap H_{2})=\overline{r}(L).$$
    It is easy to see that $H_{1}'$ and $H_{2}'$ are also factor-free. Since $H_{1}'$ and $H_{2}'$ are subgroups in a free product $G_{1}/T*G_{2}/T$ (without an amalgamated subgroup), we can use theorem 1. According to (\ref{eq2}), we have:

       \begin {equation} \label{eq5}
        \overline{r}(H_{1}'\cap H_{2}') \leq 2\frac{q_{f}^{*}}{q_{f}^{*}-2} \overline{r}(H_{1}') \overline{r}(H_{2}').
       \end {equation}

           Furthermore, $L$ is a subgroup in a free group $H_{1}'\cap H_{2}'$, so, according to Shreier formula \cite{5}, we have $$\overline{r}(L)=|H_{1}'\cap H_{2}':L|\:\overline{r}(H_{1}'\cap H_{2}')$$ Together with (\ref{eq5}), this implies the desired inequality (\ref{eq4}).     $\blacksquare$

    \begin {lemmas} \label {lemma 2}
     \begin {equation} \label{eq6}
     |H_{1}'\cap H_{2}':L| \leqslant |T|
     \end {equation}
     More exactly, all the right cosets of $L$ in $H_{1}'\cap H_{2}'$ have the form
     \begin {equation} \label{eq7}
     \varphi (H_{1}\cap H_{2}t), \quad t \in T: \quad H_{1}\cap H_{2}t \neq \varnothing
     \end{equation}
    \end {lemmas}
     $\square$    {\it \bf a)} $ \bigcup_{t \in T} \quad \varphi (H_{1} \cap H_{2}t) = H_{1}'\cap H_{2}' $
     \par The inclusion $\subseteq$ holds, since $\varphi(H_{1})=H_{1}'$, \: $\varphi (H_{2}t)=\varphi(H_{2})=H_{2}'$; \quad the inclusion $\supseteq$ holds, since
     $$g'\in H_{1}' \cap H_{2}' \quad \Rightarrow \quad g' = \varphi(h_{1})= \varphi (h_{2}), \quad h_{1} = h_{2}t, \quad h_{1} \in  H_{1}, \quad h_{2} \in H_{2}, \quad t \in T$$ $$ \Rightarrow \quad g' \in \varphi (H_{1} \cap H_{2}t) $$

     {\it \bf b)} For each $t\in T$ $\varphi (H_{1}\cap H_{2}t)$ belongs to one coset of $L$, since
     $$a',b'\in \varphi (H_{1} \cap H_{2}t)\quad \Rightarrow \quad a'=\varphi(a), \quad b'= \varphi(b), \quad a, b \in H_{1}\cap H_{2}t, \quad ab^{-1}\in H_{1} \cap H_{2}$$ $$\Rightarrow \quad a'b'^{-1}=\varphi(ab^{-1})\in \varphi (H_{1}\cap H_{2}) = L $$

     {\it \bf c)} Finally, $L\varphi(H_{1}\cap H_{2}t)= \varphi(H_{1}\cap H_{2})\varphi(H_{1}\cap H_{2}t)=\varphi (H_{1} \cap H_{2}t)$
     $\blacksquare$

     \par
       (\ref{eq6}) together with (\ref {eq4}) imply that
     $r(H_{1} \cap H_{2})$ is finite and the inequality (\ref{eq3}) holds.
    \par

    \subsection* {3. Graph of a subgroup in a free product of groups.}
     To prove that the first estimate in (\ref{eq3}) is sharp in the case mentioned above, we will use the same graph-theoretic approach as in \cite{3}, \cite{8}, \cite{9}, \cite{10}. Below we repeat the definitions from \cite{3} in the case of factor-free subgroups.

    Suppose $H\subseteq G=\prod\nolimits^{*} G_{\alpha}$ is a factor-free subgroup of $G$. We define a graph $\Psi^{*}(H)$ associated with subgroup $H$:
    \par  Vertices of $\Psi^{*}(H)$ can belong to 2 different types:
    \par 1) {\it Primary} vertices correspond to right cosets of $H$ in $G$;
    \par 2) {\it Secondary} vertices, associated with the factor $G_{\alpha}$, correspond to eqivalence classes of right cosets of $H$ in $G$  under the following equivalence relation: $Hg_{1}\sim Hg_{2}$ if $\exists c \in G_{\alpha}: Hg_{1}c = Hg_{2}$ (we will denote equivalence class as $[Hg_{1}]_{\alpha}$).
    \par  Furthermore, each primary vertex $Hg_{1}$ is connected by one edge to the secondary vertex $[Hg_{1}]_{\alpha}$ for all $\alpha$, and there are no other edges in $\Psi^{*}(H)$.

    \par  We assume that all edges of $\Psi^{*}(H)$ are oriented (if $e$ is an edge, then $e^{-1}$ denotes this edge with opposite orientation). \par Now we assign a label to every edge of $\Psi^{*}(H)$ as follows:
    \par Suppose $e_{1}$ is an edge with one end in a secondary vertex $[Hg_{1}]_{\alpha}$, then $\varphi(e_{1}) \in G_{\alpha}$ and $\varphi(e_{1}^{-1}) = \varphi(e_{1})^{-1}$; furthermore, if $Hg_{1}, Hg_{2}$ are primary vertices, connected with a secondary vertex $[Hg_{1}]_{\alpha}$ by edges $e_{1}, e_{2}$ with labels $\varphi(e_{1}), \varphi(e_{2})$ respectively, then

    \begin {equation}\label{eq8}
    Hg_{1}\varphi(e_{1})\varphi(e_{2})^{-1} = Hg_{2}
    \end {equation}

    (here we suppose that both edges are oriented from the primary vertices to the secondary vertex).

    \par
     We call {\it base vertex} primary vertex, corresponding to subgroup $H$.
          \par  {\it Label} of a path $p=e_{1}\ldots e_{m}$ in $\Psi^{*}(H)$ is $\varphi(e_{1})\ldots \varphi(e_{m}).$ \par
     A path is called {\it reduced} if it contains no subpath of the form $dd^{-1}$, where $d$ is an edge.
    \par
     Let $\Psi(H)$ be the minimal connected subgraph of $\Psi^{*}(H)$, containing all reduced cycles and the base vertex of $\Psi^{*}(H)$.
    \par
     The following 2 lemmas are also proved in \cite{3}.

     \begin {lemmas} \label {lemma 3}
      Suppose $H\subseteq G=\prod\nolimits^{*} G_{\alpha}$, $w$ is a nonempty reduced word in the alphabet $\bigcup G_{\alpha}$. Then $w \in H$ if and only if there is a reduced closed path in $\Psi(H)$ ending in the base vertex with the label $w$.
    \end {lemmas}
    \par $\square$ $(\Leftarrow)$ \: $w$ is a label of a path, which begins and ends in the primary vertex $H$, so, according to (\ref{eq8}), $Hw=H$, and that means $w\in H$. \par

    $(\Rightarrow) \: w=s_{1}\ldots s_{k},$ \: $s_{i}\in G_{\alpha_{i}}$. According to the definition of $\Psi^{*}(H)$, it contains primary vertices $H, Hs_{1}, Hs_{1}s_{2},\ldots , Hs_{1}s_{2}\ldots s_{k} = H$ and a path beginning in $Hs_{1}s_{2}\ldots s_{j-1}$ and ending in $Hs_{1}s_{2}\ldots s_{j-1}s_{j}$ with label $s_{j}$ (this path has 2 edges, both incident to the secondary vertex $[Hs_{1}s_{2}\ldots s_{j-1}]_{\alpha_{j}}=[Hs_{1}s_{2}\ldots s_{j-1}s_{j}]_{\alpha_{j}}$, it is easy to see that it has label $s_{j}$, since $H$ is factor-free). It is also easy to check that this path is in $\Psi(H)$.  $\blacksquare$

    \par
    \vspace {0.2 cm}
     Suppose $T$ is a maximal subtree in $\Psi(H)$, and $C$ is the set of edges of $\Psi(H)$, not belonging to $T$. For every edge $e$ from $C$ we can find (unique) path $q_{e}$ in $T$ beginning in the base vertex and ending in the beginning of $e$ and (unique) path $r_{e}$ in $T$ beginning in the end of $e$ and ending in the base vertex.
    %    \begin {equation}\label{eq9}
    %    q_{e}er_{e}, e\in C
    %    \end {equation}

    \begin {lemmas} \label {lemma 4}
       Suppose $H\subseteq G=\prod\nolimits^{*} G_{\alpha}$ is factor-free. Then $H$ is freely generated by
       \begin {equation}\label{eq10}
      \varphi(q_{e}er_{e}), e\in C.
       \end {equation}
       Furthermore, $H$ has finite rank if and only if the graph $\Psi(H)$ is finite, and in this case

    \begin {equation}\label{eq11}
     \overline{r}(H) = - \chi (\Psi (H)),
     \end {equation}
    where $\chi (\Psi(H))$ is the {\it Euler characteristic} of $\Psi(H)$, equal to the the number of vertices in $\Psi(H)$ minus the number of edges in $\Psi(H)$.
    \end {lemmas}

    \par $\square$ The full proof of this lemma can be found in \cite{3}; it uses the fact that the fundamental group of $\Psi(H)$ is freely generated by the paths $q_{e}er_{e}, e\in C.$
     Furthermore, $$ - \chi (\Psi (H)) = |C| - 1, $$ since all the vertices of $\Psi (H)$ are in $T$ and the number of edges in $T$ is less than the number of vertices in $T$ by 1 ($T$ is a tree). According to (\ref{eq10}), $|C|=r(H)$ and (\ref{eq11}) holds.\quad $\blacksquare$

    \par
     It is easy to see that, if $H$ has finite index in $G$, then $\Psi(H)=\Psi^{*}(H).$ \par
    It is also easy to see that, if $H$ is a normal subgroup of $G$, then $\Psi^{*}(H)$ corresponds to the Cayley graph of $G/H$, with primary vertices of $\Psi^{*}(H)$ corresponding to the vertices of the Cayley graph of $G/H$.

    \subsection* {4. Proof of unimprovability of the estimate.}
    Consider again the factorization  $$\varphi : G = G_{1}*_{T}G_{2} \rightarrow G_{1}/T*G_{2}/T = G_{1}'*G_{2}' = G'.$$
    It is given that $G'$ contains an involution and $G' \ncong \mathbb{Z}_{2}*\mathbb{Z}_{2}$.
    The following lemma is also proved in \cite{4}:
    \begin {lemmas} \label {lemma 5}
     Suppose $G'=G_{1}'*G_{2}'$, \quad $q_{f}^{*}$ is the minimum of orders $> 2$ of subgroups of groups $G_{1}'$, $G_{2}'$. Suppose also $G'$ contains an involution and $G' \ncong \mathbb{Z}_{2}*\mathbb{Z}_{2}$. Then either $q_{f}^{*} = \infty$, or $q_{f}^{*}$ is prime, or $q_{f}^{*} = 4$, and $G'$ has one of the following subgroups:
    \begin {equation} \label{eq12a}
    \mathbb{Z}_{2}*\mathbb{Z}_{2}*\mathbb{Z}_{2}, \quad \text{if} \quad q_{f}^{*} = \infty
    \end {equation}
    \begin {equation} \label{eq12b}
    \mathbb{Z}_{2}*\mathbb{Z}_{p}, \quad \text{if} \quad q_{f}^{*} = p, \quad \text{p is prime}, \quad p > 2
    \end {equation}
    \begin {equation} \label{eq12c}
    \mathbb{Z}_{2}*\mathbb{Z}_{4}  \quad \text{or} \quad  \mathbb{Z}_{2}* (\mathbb{Z}_{2} \times \mathbb{Z}_{2}), \quad  \text{if} \quad q_{f}^{*} = 4.
    \end {equation}
    \end {lemmas}
    \par
    $\square$  The first statement follows from the definition of $q_{f}^{*}$ and Sylow theorems. \par
    Furthermore, suppose $a$ is an involution in $G'$, then we can assume that it belongs to one of the factors, for instance $a \in G_{1}'$. If $q_{f}^{*} < \infty$, then let $Q$ be the subgroup of $G_{1}'$ or $G_{2}'$ with $q_{f}^{*}$ elements. Then the subgroup $\langle bab^{-1} \rangle * Q$, where $b \in G_{2}', \: b \neq 1$, has the desired form (\ref{eq12b}) or (\ref{eq12c}). \par If $q_{f}^{*} = \infty$, then the subgroup $ \langle a  \rangle * \langle bab^{-1} \rangle *  \langle cac^{-1} \rangle $, where $b,c \in G_{2}' , \: b, c\neq 1, \: b \neq c$, has the desired form (\ref{eq12a}). $\blacksquare$

    \par \vspace{0.2 cm}
    We want to prove that in the case under consideration the first estimate in (\ref{eq3}) cannot be further improved. Evidently, we can consider a subgroup of the form (\ref{eq12a}), (\ref{eq12b}) or (\ref{eq12c}) in $G'$ and its full inverse image under the homomorphism $\varphi$ in $G$ and restrict $\varphi$ on this inverse image. So it will be enough to prove that the first estimate in (\ref{eq3}) is sharp (and that means to give examples of appropriate subgroups $H_{1}$ and $H_{2}$) in the case when $G'$ is one of the groups (\ref{eq12a}), (\ref{eq12b}) or (\ref{eq12c}).  \par

    It is easy to see from the proofs of lemmas 1 and 2 that the first inequality in (\ref{eq3}) turns into equality  if and only if both inequalities (\ref{eq5}) and (\ref{eq6}) turn into equalities. We do the following: for every $n = |T|$ for each of the groups (\ref{eq12a}), (\ref{eq12b}) and (\ref{eq12c}) we construct in them (free) subgroups of finite index $H_{1}'$ and $H_{2}'$, such that (\ref{eq5}) turns into equality for both subgroups; afterwards we choose the inverse images of the free generators of these subgroups under the homomorphism $\varphi$, so that (\ref{eq6}) turns into equality for both inverse images $H_{1}$ and $H_{2}$. More formally this process is described below. \par

    The following lemma is also proved in \cite{4}.
    \begin {lemmas} \label {lemma 6}
    Suppose $H_{1}'$ and $H_{2}'$ are factor-free subgroups of finite index in $G'$, and $G'$ is one of the groups (\ref{eq12a}), (\ref{eq12b}) or (\ref{eq12c}). Then (\ref{eq5}) turns into equality if and only if
    \begin {equation} \label{eq12.5}
    |G': H_{1}' \cap H_{2}'| = |G':H_{1}'|\cdot|G':H_{2}'|
    \end {equation}
    \end {lemmas}
    $\square$
     First notice that, according to (\ref{eq11}),
     \begin{equation} \label{eq12.7}
       \overline {r}(H')= - \chi (\Psi (H')) = \frac {1}{2} \sum\limits_{v\in V(\Psi(H'))}(\mbox{deg\:} v - 2),
     \end{equation}
     where $H' \subseteq G', \quad V(\Psi(H'))$ is the set of vertices of graph $\Psi(H')$.
     Let $U(\Psi(H'))$ denote the set of primary vertices of $\Psi(H')$, and $V_{p}(\Psi(H'))$ denote the set of vertices of degree $p$ in $\Psi(H')$. Obviously, $|U(\Psi(H'))| = |G':H'|$. \par
     Below $H'$ denotes one of the subgroups $H_{1}'$, $H_{2}'$ or $H_{1}' \cap H_{2}'$ in $G'$ (all of them are factor-free and have finite index in $G'$). \par
     First consider the case when $G'$ has the form (\ref{eq12b}) or (\ref{eq12c}). Then the primary vertices of $\Psi(H')$ have degree 2, the secondary vertices, corresponding to the factor $\mathbb{Z}_{2}$, also have degree 2, and the secondary vertices, corresponding to the other factor, have degree $p$, where $p$ is $4$ or prime. Notice that the vertices of degree 2 do not influence on the sum on the right-hand side of (\ref{eq12.7}), so we obtain
     $$ \overline {r}(H')=  \frac {p - 2}{2} |V_{p}(\Psi(H'))|.$$
     Furthermore, each primary vertex is connected with exactly one secondary vertex of degree $p$, therefore, $$|G':H'| = |U(\Psi(H'))| = p|V_{p}(\Psi(H'))| = \frac {2p}{p-2}\overline {r}(H').$$
     Substituting all the 3 indexes in (\ref{eq12.5}) according to the last formula and considering that $q^{*} = p$, we obtain that in this case (\ref{eq12.5}) is equivalent to the equality in (\ref{eq5}).
     \par
     Now consider the remaining case when $G' = \mathbb{Z}_{2}*\mathbb{Z}_{2}*\mathbb{Z}_{2}$, $q^{*} = \infty$. In this case all the secondary vertices of $\Psi(H')$ have degree 2, and all the primary vertices have degree 3, so we obtain:
     $$ \overline {r}(H') = \frac{|V_{3}(\Psi(H'))|}{2} =   \frac{|U(\Psi(H'))|}{2} =  \frac{|G':H'|}{2},$$
      therefore, $|G':H'| = 2 \overline {r}(H')$. Substitutung all the 3 indexes in (\ref{eq12.5}) according to the last formula and considering that $\frac{q^{*}}{q^{*}-2} = 1$, we obtain that in this case also (\ref{eq12.5}) is equivalent to the equality in (\ref{eq5}). $\blacksquare$ \par
    The following lemma is well-known:
    \begin {lemmas} \label {lemma 7}
    Suppose
    \begin {equation} \label{eq13}
    H_{1}' \triangleleft G', \quad H_{1}'H_{2}' = G'
    \end {equation}
    Then $|G': H_{1}' \cap H_{2}'| = |G':H_{1}'|\cdot|G':H_{2}'|$
    \end {lemmas}
    $\square $
    $$ \quad G'/H_{1}' = (H_{1}'H_{2}')/H_{1}' \cong H_{2}'/(H_{1}'\cap H_{2}'),$$ so
    $|G':H_{1}'| = |H_{2}':H_{1}' \cap H_{2}'|$, therefore, $$|G': H_{1}' \cap H_{2}'| = |G':H_{2}'|\cdot |H_{2}':H_{1}' \cap H_{2}'| = |G':H_{2}'|\cdot |G':H_{1}'| \quad \blacksquare$$

    We obtain from lemmas 6 and 7 that in the case under consideration the condition (\ref{eq13}) is sufficient for equality in (\ref{eq5}).

    Consider now the inequality (\ref{eq6}). It follows from lemma 2 that it turns into equality if and only if
     \begin {equation} \label{eq14}
    H_{1}\cap H_{2}t \neq \varnothing \quad \forall t\in T
    \end {equation}
    And (\ref{eq14}) holds when
    \begin {equation} \label{eq15}
    T\subseteq H_{1}H_{2}.
    \end {equation}
    Therefore, in the examples constructed below it is enough to prove that (\ref{eq13}) and (\ref{eq15}) hold, and that both subgroups $H_{1}'$ and $H_{2}'$ are factor-free and have finite index in $G'$, then the first inequality in (\ref{eq3}) turns into equality.

    \par

    Now we proceed to constructing the desired examples. Let $|T|=n$.
    \par \textbf{Case 1.}
    \quad Suppose first that $$G' = \mathbb{Z}_{2}* \mathbb{Z}_{p} \cong \langle a, b \:|\: a^{p}=b^{2}=1 \rangle,$$ where $p$ is prime, $p > 2$. \par
    First we construct subgroup $H_{1}'$. \par
     Consider $$G_{0}' = \langle \langle (ab)^{6} \rangle \rangle \subseteq G', \quad R = G'/ G_{0}' \cong \langle a, b \:|\: a^{p}=b^{2}=(ab)^6=1\rangle = T(p,2,6). $$ $R$ is a triangle group; it is well-known that it is infinite, since $$\frac{1}{p}+\frac{1}{2}+\frac{1}{6} \leq 1 \quad \text {when} \: \: p \geq 3. $$ It is also well-known that triangle groups are residually finite (\cite{6}). Therefore, for each finite subset $M$ of group $R$ there exists a homomorphism from $R$ on a finite group, injective on the set $M$. \par

    Consider the following elements of the group $G_{0}'$:
      \begin{equation}\label{eq16}
      w_{1}' = (ab)^{6}, w_{2}' = ((ab)^{6})^{ba^{-1}ba^{-2}}, w_{3}' = ((ab)^{6})^{(ba^{-1}ba^{-2})^{2}}, ..., w_{n}' = ((ab)^{6})^{(ba^{-1}ba^{-2})^{n-1}}
      \end{equation}
       It is easy to see that all these elements are different. \par
       Consider part $I$ of the graph $\Psi(G_{0}')$, which consists of all secondary vertices of this graph, which lie on the reduced paths with labels (\ref{eq16}), together with the incident edges and their other ends --- primary vertices. Let $M$ be the set of the right cosets, corresponding to the primary vertices of $I$, then obviously $M$ is the subset of quotient group $R = G'/ G_{0}'$ and $1, a, b \in M$. \par
       Part $I$ of the graph $\Psi(G_{0}')$ in the case $n = 3$ is represented on picture 1.
       \par

         %\footnotesize
         \scriptsize

\vspace{0.5 cm}
\hspace{0.7 cm}
\xy
0;/r.17 pc/:
(0,0)*{\circlearrowright}; (10,-10)*{\circ} **\dir{-};
(10,-10)*{\circ}; (20,-20)*{\circlearrowright} **\dir{-};
(5,-5)*{\bullet};
(15,-15)*{\bullet};
(20,-20)*{\circlearrowright}; (20,-30)*{\circ} **\dir{-};
(20,-30)*{\circ}; (20,-40)*{\circlearrowright} **\dir{-};
(20,-25)*{\bullet};
(20,-35)*{\bullet};
(20,-40)*{\circlearrowright}; (10,-50)*{\circ} **\dir{-};
(10,-50)*{\circ}; (0,-60)*{\circlearrowright} **\dir{-};
(15,-45)*{\bullet};
(5,-55)*{\bullet};
(0,0)*{\circlearrowright}; (-10,-10)*{\circ} **\dir{-};
(-10,-10)*{\circ}; (-20,-20)*{\circlearrowright} **\dir{-};
(-5,-5)*{\bullet};
(-15,-15)*{\bullet};
(-13,-17)*{\ast};
(-20,-20)*{\circlearrowright}; (-20,-30)*{\circ} **\dir{-};
(-20,-30)*{\circ}; (-20,-40)*{\circlearrowright} **\dir{-};
(-20,-25)*{\bullet};
(-20,-35)*{\bullet};
(-20,-40)*{\circlearrowright}; (-10,-50)*{\circ} **\dir{-};
(-10,-50)*{\circ}; (0,-60)*{\circlearrowright} **\dir{-};
(-15,-45)*{\bullet};
(-5,-55)*{\bullet};
(0,0)*{\circlearrowright}; (0,5)*{\bullet} **\dir{-};
(-20,-20)*{\circlearrowright}; (-25,-15)*{\bullet} **\dir{-};
(-20,-40)*{\circlearrowright}; (-25,-45)*{\bullet} **\dir{-};
(0,-60)*{\circlearrowright}; (0,-65)*{\bullet} **\dir{-};
(20,-40)*{\circlearrowright}; (30,-50)*{\circ} **\dir{-};
(30,-50)*{\circ}; (40,-60)*{\circlearrowright} **\dir{-};
(25,-45)*{\bullet};
(35,-55)*{\bullet};
(40,0)*{\circlearrowright}; (50,-10)*{\circ} **\dir{-};
(50,-10)*{\circ}; (60,-20)*{\circlearrowright} **\dir{-};
(45,-5)*{\bullet};
(55,-15)*{\bullet};
(60,-20)*{\circlearrowright}; (60,-30)*{\circ} **\dir{-};
(60,-30)*{\circ}; (60,-40)*{\circlearrowright} **\dir{-};
(60,-25)*{\bullet};
(60,-35)*{\bullet};
(60,-40)*{\circlearrowright}; (50,-50)*{\circ} **\dir{-};
(50,-50)*{\circ}; (40,-60)*{\circlearrowright} **\dir{-};
(55,-45)*{\bullet};
(45,-55)*{\bullet};
(40,0)*{\circlearrowright}; (30,-10)*{\circ} **\dir{-};
(30,-10)*{\circ}; (20,-20)*{\circlearrowright} **\dir{-};
(35,-5)*{\bullet};
(25,-15)*{\bullet};
(40,0)*{\circlearrowright}; (40,5)*{\bullet} **\dir{-};
(40,-60)*{\circlearrowright}; (40,-65)*{\bullet} **\dir{-};
(60,-20)*{\circlearrowright}; (65,-15) **\dir{-};
(60,-40)*{\circlearrowright}; (65,-45) **\dir{-};
(70,-30)*{\ldots};
(100,0)*{\circlearrowright}; (110,-10)*{\circ} **\dir{-};
(110,-10)*{\circ}; (120,-20)*{\circlearrowright} **\dir{-};
(105,-5)*{\bullet};
(115,-15)*{\bullet};
(120,-20)*{\circlearrowright}; (120,-30)*{\circ} **\dir{-};
(120,-30)*{\circ}; (120,-40)*{\circlearrowright} **\dir{-};
(120,-25)*{\bullet};
(120,-35)*{\bullet};
(120,-40)*{\circlearrowright}; (110,-50)*{\circ} **\dir{-};
(110,-50)*{\circ}; (100,-60)*{\circlearrowright} **\dir{-};
(115,-45)*{\bullet};
(105,-55)*{\bullet};
(100,0)*{\circlearrowright}; (90,-10)*{\circ} **\dir{-};
(90,-10)*{\circ}; (80,-20)*{\circlearrowright} **\dir{-};
(95,-5)*{\bullet};
(85,-15)*{\bullet};
(80,-20)*{\circlearrowright}; (80,-30)*{\circ} **\dir{-};
(80,-30)*{\circ}; (80,-40)*{\circlearrowright} **\dir{-};
(80,-25)*{\bullet};
(80,-35)*{\bullet};
(80,-40)*{\circlearrowright}; (90,-50)*{\circ} **\dir{-};
(90,-50)*{\circ}; (100,-60)*{\circlearrowright} **\dir{-};
(85,-45)*{\bullet};
(95,-55)*{\bullet};
(100,0)*{\circlearrowright}; (100,5)*{\bullet} **\dir{-};
(120,-20)*{\circlearrowright}; (125,-15)*{\bullet} **\dir{-};
(120,-40)*{\circlearrowright}; (125,-45)*{\bullet} **\dir{-};
(100,-60)*{\circlearrowright}; (100,-65)*{\bullet} **\dir{-};
(80,-20)*{\circlearrowright}; (75,-15) **\dir{-};
(80,-40)*{\circlearrowright}; (75,-45) **\dir{-};
(0,-30)*{1};
(40,-30)*{2};
(100,-30)*{n};
(17,-13)*{{\bf u}};
(23,-13)*{{\bf v}};
(33,-3)*{{\bf w}};
(17,-26)*{{\bf y}};
(37,2)*{a};
(43,2)*{a};
(40,-3)*{a};
(97,2)*{a};
(103,2)*{a};
(100,-3)*{a};
(-3,2)*{a};
(3,2)*{a};
(0,-3)*{a};
(-23,-38)*{a};
(-17,-38)*{a};
(-20,-43)*{a};
(23,-38)*{a};
(17,-38)*{a};
(20,-43)*{a};
(63,-38)*{a};
(57,-38)*{a};
(60,-43)*{a};
(83,-38)*{a};
(77,-38)*{a};
(80,-43)*{a};
(123,-38)*{a};
(117,-38)*{a};
(120,-43)*{a};
(-20,-17)*{a};
(-23,-22)*{a};
(-17,-22)*{a};
(20,-17)*{a};
(23,-22)*{a};
(17,-22)*{a};
(60,-17)*{a};
(63,-22)*{a};
(57,-22)*{a};
(80,-17)*{a};
(83,-22)*{a};
(77,-22)*{a};
(120,-17)*{a};
(123,-22)*{a};
(117,-22)*{a};
(0,-57)*{a};
(-3,-62)*{a};
(3,-62)*{a};
(40,-57)*{a};
(37,-62)*{a};
(43,-62)*{a};
(100,-57)*{a};
(97,-62)*{a};
(103,-62)*{a};
\endxy
\vspace{0.1 cm}
\par
\begin{center}
Picture 1.
\end{center}
\par
\normalsize
\vspace{0.1 cm}

    (On this picture primary vertices are represented as $\bullet$, secondary vertices, corresponding to the factor $\langle b  \rangle _{2}$ --- as $\circ$,
    and secondary vertices, corresponding to the factor $ \langle a  \rangle _{3}$ --- as $\circlearrowright$. The label of a nontrivial path with 2 edges, both incident to the secondary vertex, corresponding to the factor $\langle b  \rangle _{2}$, is equal to $b$ (and not represented on the picture). The label of a nontrivial path with 2 edges, both incident to the secondary vertex, corresponding to the factor $ \langle a  \rangle _{3}$, is equal to $a$ or $a^{2}$, depending on the labels on the picture next to this secondary vertex and on the direction we go around this secondary vertex (in all our examples we always go around the vertices clockwise). For example, the label of the path (with 2 edges) from ${\bf v}$ to ${\bf w}$ is equal to $b$; the label of the path from ${\bf u}$ to ${\bf v}$ is equal to $a$, and the label of the path from ${\bf u}$ to ${\bf y}$ is equal to $a^{2}$. The base vertex is marked with a symbol $\ast$ next to it.)

    According to the facts mentioned above, there exists a surjective homomorphism $\pi: R \rightarrow S = R/R_{0}$, injective on the set $M$, where $S$ is a finite group. Then $S \cong G'/H_{1}'$, $H_{1}' \lhd G'$, $G_{0}' \subseteq H_{1}'$. Subgroup $H_{1}'$ has finite index in $G'$, since $S$ is finite. Moreover, $H_{1}'$ is factor-free (and therefore, $H_{1}'$ is free), since otherwise normal subgroup $H_{1}'$ would contain $a$ or $b$, but it is impossible, since the elements $1, a, b \in M$ and are therefore injectively mapped on $S$. \par

        It is easy to see that by choosing the maximal subtree in the part $I$ of the graph $\Psi (G_{0}')$ in the appropriate way we can make all $n$ elements (\ref{eq16}) belong to the free generators of $G_{0}'$.
       Moreover, we can make all elements (\ref{eq16}) belong to the free generators of $H_{1}'$. It follows from the injectivity of $\pi$ on $M$: part of $\Psi(H_{1}')$, corresponding to the elements (\ref{eq16}), will be the same as part $I$ of $\Psi(G_{0}')$, and by choosing the maximal subtree in the same way in $\Psi (H_{1}')$ we obtain the desired condition. Thus, we can suppose that all elements (\ref{eq16}) belong to the free generators of $H_{1}'$.
      Notice also that $(ba)^{6}=((ab)^6)^{b}\in H_{1}'.$

        Now we construct subgroup $H_{2}'$. \par
         Consider subgroup $K\subseteq G'$, (freely) generated by the following elements:
         \begin{equation}\label{eq18}
         w_{1}', ..., w_{n}' \quad \text{from} \: (\ref{eq16}), \quad (ba)^{5}, \quad (ba)^{2}(ba^{-1})^{5}(ba)^{2}.
         \end{equation}
          The graph $\Psi(K)$ in the case $n = 3$ is represented on picture 2.
         \par

%\footnotesize
\scriptsize
\vspace{0.5 cm}
\hspace{0.5 cm}
\xy
0;/r.15 pc/:
(0,0)*{\circlearrowright}; (-10,10)*{\circ} **\dir{-};
(-10,10)*{\circ}; (-20,20)*{\circlearrowright} **\dir{-};
(-22,16)*{a};
(-5,5)*{\bullet};
(-15,15)*{\bullet};
(-13,-17)*{\ast};
(-20,20)*{\circlearrowright}; (-30,20)*{\circ} **\dir{-};
(-30,20)*{\circ}; (-40,20)*{\circlearrowright} **\dir{-};
(-25,20)*{\bullet};
(-35,20)*{\bullet};
(-20,-20)*{\circlearrowright}; (-30,-10)*{\circ} **\dir{-};
(-30,-10)*{\circ}; (-40,0)*{\circlearrowright} **\dir{-};
(-37,2)*{a};
(-25,-15)*{\bullet};
(-35,-5)*{\bullet};
(-40,0)*{\circlearrowright}; (-40,10)*{\circ} **\dir{-};
(-40,10)*{\circ}; (-40,20)*{\circlearrowright} **\dir{-};
(-40,5)*{\bullet};
(-40,15)*{\bullet};
(-40,20)*{\circlearrowright}; (-50,30)*{\circ} **\dir{-};
(-50,30)*{\circ}; (-60,40)*{\circlearrowright} **\dir{-};
(-45,25)*{\bullet};
(-55,35)*{\bullet};
(-60,40)*{\circlearrowright}; (-70,40)*{\circ} **\dir{-};
(-70,40)*{\circ}; (-70,50)*{\circlearrowright} **\dir{-};
(-73,53)*{a};
(-65,40)*{\bullet};
(-70,45)*{\bullet};
(-60,40)*{\circlearrowright}; (-60,50)*{\circ} **\dir{-};
(-60,50)*{\circ}; (-70,50)*{\circlearrowright} **\dir{-};
(-60,45)*{\bullet};
(-65,50)*{\bullet};
%
%(-70,50)*{\circlearrowright}; (-65,45) **\dir{-};
%
%
%
(0,0)*{\circlearrowright}; (10,-10)*{\circ} **\dir{-};
(10,-10)*{\circ}; (20,-20)*{\circlearrowright} **\dir{-};
(5,-5)*{\bullet};
(15,-15)*{\bullet};
(20,-20)*{\circlearrowright}; (20,-30)*{\circ} **\dir{-};
(20,-30)*{\circ}; (20,-40)*{\circlearrowright} **\dir{-};
(20,-25)*{\bullet};
(20,-35)*{\bullet};
(20,-40)*{\circlearrowright}; (10,-50)*{\circ} **\dir{-};
(10,-50)*{\circ}; (0,-60)*{\circlearrowright} **\dir{-};
(15,-45)*{\bullet};
(5,-55)*{\bullet};
(0,0)*{\circlearrowright}; (-10,-10)*{\circ} **\dir{-};
(-10,-10)*{\circ}; (-20,-20)*{\circlearrowright} **\dir{-};
(-5,-5)*{\bullet};
(-15,-15)*{\bullet};
(-20,-20)*{\circlearrowright}; (-20,-30)*{\circ} **\dir{-};
(-20,-30)*{\circ}; (-20,-40)*{\circlearrowright} **\dir{-};
(-17,-38)*{a};
(-20,-25)*{\bullet};
(-20,-35)*{\bullet};
(-20,-40)*{\circlearrowright}; (-10,-50)*{\circ} **\dir{-};
(-10,-50)*{\circ}; (0,-60)*{\circlearrowright} **\dir{-};
(0,-56)*{a};
(-15,-45)*{\bullet};
(-5,-55)*{\bullet};
%
%(0,0)*{\circlearrowright}; (0,5)*{\bullet} **\dir{-};
%(-20,-20)*{\circlearrowright}; (-25,-15)*{\bullet} **\dir{-};
%(-20,-40)*{\circlearrowright}; (-25,-45)*{\bullet} **\dir{-};
%(0,-60)*{\circlearrowright}; (0,-65)*{\bullet} **\dir{-};
%
%
(20,-40)*{\circlearrowright}; (30,-50)*{\circ} **\dir{-};
(30,-50)*{\circ}; (40,-60)*{\circlearrowright} **\dir{-};
(40,-56)*{a};
(25,-45)*{\bullet};
(35,-55)*{\bullet};
(40,-4)*{a};
(40,0)*{\circlearrowright}; (50,-10)*{\circ} **\dir{-};
(50,-10)*{\circ}; (60,-20)*{\circlearrowright} **\dir{-};
(45,-5)*{\bullet};
(55,-15)*{\bullet};
(60,-20)*{\circlearrowright}; (60,-30)*{\circ} **\dir{-};
(60,-30)*{\circ}; (60,-40)*{\circlearrowright} **\dir{-};
(60,-25)*{\bullet};
(60,-35)*{\bullet};
(60,-40)*{\circlearrowright}; (50,-50)*{\circ} **\dir{-};
(50,-50)*{\circ}; (40,-60)*{\circlearrowright} **\dir{-};
(55,-45)*{\bullet};
(45,-55)*{\bullet};
(40,0)*{\circlearrowright}; (30,-10)*{\circ} **\dir{-};
(30,-10)*{\circ}; (20,-20)*{\circlearrowright} **\dir{-};
(35,-5)*{\bullet};
(25,-15)*{\bullet};
%
%(40,0)*{\circlearrowright}; (40,5)*{\bullet} **\dir{-};
%(40,-60)*{\circlearrowright}; (40,-65)*{\bullet} **\dir{-};
%
(60,-20)*{\circlearrowright}; (65,-15) **\dir{-};
(60,-40)*{\circlearrowright}; (65,-45) **\dir{-};
(70,-30)*{\ldots};
(100,0)*{\circlearrowright}; (110,-10)*{\circ} **\dir{-};
(110,-10)*{\circ}; (120,-20)*{\circlearrowright} **\dir{-};
(100,-4)*{a};
(117,-22)*{a};
(105,-5)*{\bullet};
(115,-15)*{\bullet};
(120,-20)*{\circlearrowright}; (120,-30)*{\circ} **\dir{-};
(120,-30)*{\circ}; (120,-40)*{\circlearrowright} **\dir{-};
(117,-38)*{a};
(120,-25)*{\bullet};
(120,-35)*{\bullet};
(120,-40)*{\circlearrowright}; (110,-50)*{\circ} **\dir{-};
(110,-50)*{\circ}; (100,-60)*{\circlearrowright} **\dir{-};
(100,-56)*{a};
(115,-45)*{\bullet};
(105,-55)*{\bullet};
(100,0)*{\circlearrowright}; (90,-10)*{\circ} **\dir{-};
(90,-10)*{\circ}; (80,-20)*{\circlearrowright} **\dir{-};
(95,-5)*{\bullet};
(85,-15)*{\bullet};
(80,-20)*{\circlearrowright}; (80,-30)*{\circ} **\dir{-};
(80,-30)*{\circ}; (80,-40)*{\circlearrowright} **\dir{-};
(80,-25)*{\bullet};
(80,-35)*{\bullet};
(80,-40)*{\circlearrowright}; (90,-50)*{\circ} **\dir{-};
(90,-50)*{\circ}; (100,-60)*{\circlearrowright} **\dir{-};
(85,-45)*{\bullet};
(95,-55)*{\bullet};
%
%(100,0)*{\circlearrowright}; (100,5)*{\bullet} **\dir{-};
%(120,-20)*{\circlearrowright}; (125,-15)*{\bullet} **\dir{-};
%(120,-40)*{\circlearrowright}; (125,-45)*{\bullet} **\dir{-};
%(100,-60)*{\circlearrowright}; (100,-65)*{\bullet} **\dir{-};
%
(80,-20)*{\circlearrowright}; (75,-15) **\dir{-};
(80,-40)*{\circlearrowright}; (75,-45) **\dir{-};
(0,-30)*{1};
(40,-30)*{2};
(100,-30)*{n};
(-37,-7)*{{\bf u}};
(-43,5)*{{\bf v}};
(-4,0)*{a};
(3,2)*{a};
(0,-4)*{a};
(23,-38)*{a};
(17,-38)*{a};
(20,-43)*{a};
(63,-38)*{a};
(57,-38)*{a};
(60,-43)*{a};
(83,-38)*{a};
(77,-38)*{a};
(80,-43)*{a};
(-20,-16)*{a};
(-23,-22)*{a};
(-17,-22)*{a};
(20,-16)*{a};
(23,-22)*{a};
(17,-22)*{a};
(60,-16)*{a};
(63,-22)*{a};
(57,-22)*{a};
(80,-16)*{a};
(83,-22)*{a};
(77,-22)*{a};
(0,-64)*{a^{-1}};
(40,-64)*{a^{-1}};
(100,-64)*{a^{-1}};
(40,4)*{a^{-1}};
(100,4)*{a^{-1}};
(125,-16)*{a^{-1}};
(124,-43)*{a^{-1}};
(-24,-44)*{a^{-1}};
(-43,18)*{a};
(-37,17)*{a};
(-37,23)*{a};
(-44,-4)*{a^{-1}};
(-16,24)*{a^{-1}};
(-62,42)*{a};
(-62,37)*{a};
(-57,42)*{a};
(-65,47)*{a^{-1}};
\endxy
\vspace{0.1 cm}
\begin{center}
Picture 2.
\end{center}
\vspace{0.1 cm}
\normalsize

            (Here the notations are the same as on picture 1; for example, the label of the path (with 2 edges) from ${\bf u}$ to ${\bf v}$ is equal to $a^{-1}$.)

            \par
           $K$ is a factor-free subgroup of $G'$, but $K$ has infinite index in $G'$. We can construct a subgroup $H_{2}'\subseteq G'$, such that $K\subseteq H_{2}'$, $H_{2}'$ is also factor-free in $G'$ and $H_{2}'$ has finite index in $G'$. To do this suppose $J$ is the set of those primary vertices of $\Psi^{*}(K)$, which do not belong to $\Psi(K)$, but are connected by an edge with a secondary vertex of $\Psi(K)$. Suppose $Z=\{z_{j}, j=1 ... |J|\}$ is the set of paths lying in a (fixed) maximal subtree of $\Psi^{*}(K)$ beginning in the base vertex of $\Psi^{*}(K)$ and ending in a vertex of $J$, one path for each vertex of $J$.
            Let
        \begin{equation}\label{eq17}
         s_{1} = baba^{-2}b, \: s_{2} = ba^{3}ba^{-4}b, ..., \: s_{(p-1)/2} = ba^{p-2}ba^{-(p-1)}b
        \end{equation}
            Then as the generators of subgroup $H_{2}'$ we take the union of the generators (\ref{eq18}) of $K$ and the elements
          $$s_{i}^{z_{j}}, \quad i = 1 ... (p-1)/2, \: j=1 ... J. $$
          It is easy to see that $\Psi^{*}(H_{2}')=\Psi(H_{2}')$ and this graph has a finite number of primary vertices, therefore, $H_{2}'$ has finite index in $G'$; other conditions mentioned above are also obviously satisfied.

     Now we show that $H_{1}'H_{2}' = G'.$ Indeed, $$(ba)^{6} \in H_{1}', \: (ba)^{5} \in H_{2}' \Rightarrow (ba)^{6}, \: (ba)^{5} \in H_{1}'H_{2}' \Rightarrow ba \in H_{1}'H_{2}'.$$ Furthermore, $$(ba)^{2}(ba^{-1})^{5}(ba)^{2} \in H_{2}' \Rightarrow (ba)^{2}(ba^{-1})^{5}(ba)^{2}, \: ba \in H_{1}'H_{2}' \Rightarrow (ba^{-1})^{5} \in H_{1}'H_{2}'.$$ Moreover, $$(ba^{-1})^{6} = ((ab)^{6})^{-1} \in H_{1}' \Rightarrow (ba^{-1})^{5}, \: (ba^{-1})^{6} \in H_{1}'H_{2}' \Rightarrow ba^{-1} \in H_{1}'H_{2}'.$$ Thus, $a^{2} \in H_{1}'H_{2}'$, but $a^{p} = 1,$ $p$ is odd, therefore, $a, b \in H_{1}'H_{2}'$, so $H_{1}'H_{2}' = G'.$ \par

    Now we choose the inverse images of free generators of subgroups $H_{1}'$ and $H_{2}'$ under homomorphism $\varphi$ as following. We choose arbitrary inverse images $w_{1},...,w_{n} \in G$ of the elements $w_{1}',...,w_{n}'$ (from (\ref{eq16})) considered as the generators of subgroup $H_{1}'$, and we choose inverse images $w_{1}t_{1},...,w_{n}t_{n}$, where $T = \{t_{1},...,t_{n}\}$, for the same elements $w_{1}',...,w_{n}'$ considered as the generators of subgroup $H_{2}'$. We choose arbitrary inverse images of all other free generators of subgroups $H_{1}'$ and $H_{2}'$. We call $H_{1}$ and $H_{2}$ the obtained inverse images of groups $H_{1}'$ and $H_{2}'$ respectively. It is clear that (\ref{eq15}) holds. Moreover, it is proved above that (\ref{eq13}) also holds, and both subgroups $H_{1}'$ and $H_{2}'$ are factor-free and have finite index in $G'$. Thus, in this case the unimprovability of (\ref{eq3}) is proved.

    \vspace {0.2 cm}
    \par \textbf{Case 2.}  Suppose that \par
    $$ G' = \mathbb{Z}_{2}* \mathbb{Z}_{4} \cong \langle a, b \:|\: a^{4}=b^{2}=1 \rangle $$
    \par

    Subgroup $H_{1}'$ is constructed in the same way as in the previous case (the triangle group $T(4,2,6)$ is also infinite and residually finite). \par
    Now we construct subgroup $H_{2}'$. \par
    Consider subgroup $K \subseteq G'$, (freely) generated by the following elements:
    \begin {equation}\label{eq19}
    w_{1}', ..., w_{n}' \quad \text{from} \: (\ref{eq16}), \quad (ba)^{5},
     \quad (ba)^{2}(ba^{2})^{5}(ba)^{2}, \quad (ba^{2})^{2}
    \end{equation}
    The graph $\Psi(K)$ is represented below on picture 3. \par

    \par
    As in the previous case, subgroup $K$ is factor-free in $G'$,
    but $K$ has an infinite index in $G'$. Again we can construct a subgroup $H_{2}'\subseteq G'$,
    such that $K\subseteq H_{2}'$, $H_{2}'$ is also factor-free in $G'$ and
    $H_{2}'$ has finite index in $G'$. To do this suppose $J$ is the set of those primary vertices of $\Psi^{*}(K)$, which do not belong to $\Psi(K)$, but are connected by an edge with a secondary vertex of $\Psi(K)$. It is easy to see that $|J|$ is even. Suppose $Z=\{z_{j}, j=1 ... |J|\}$ is the set of paths lying in a (fixed) maximal subtree of $\Psi^{*}(K)$ beginning in the base vertex of $\Psi^{*}(K)$ and ending in a vertex of $J$, one path for each vertex of $J$.
    Then as the generators of subgroup $H_{2}'$ we take the union of the generators (\ref{eq19}) of $K$ and the elements
    $$z_{2j-1}bz_{2j}^{-1}, j= 1 ... |J|/2. $$
         It is easy to see that $\Psi^{*}(H_{2}')=\Psi(H_{2}')$ and this graph has a finite number of primary vertices, therefore, $H_{2}'$ has finite index in $G'$; other conditions mentioned above are also obviously satisfied.

    \par

    \scriptsize
    \vspace{0.5 cm}
\hspace{0.5 cm}
\xy
0;/r.15 pc/:
(0,0)*{\circlearrowright}; (-10,10)*{\circ} **\dir{-};
(-10,10)*{\circ}; (-20,20)*{\circlearrowright} **\dir{-};
(-22,16)*{a};
(-5,5)*{\bullet};
(-15,15)*{\bullet};
(-13,-17)*{\ast};
(-20,20)*{\circlearrowright}; (-30,20)*{\circ} **\dir{-};
(-30,20)*{\circ}; (-40,20)*{\circlearrowright} **\dir{-};
(-25,20)*{\bullet};
(-35,20)*{\bullet};
(-20,-20)*{\circlearrowright}; (-30,-10)*{\circ} **\dir{-};
(-30,-10)*{\circ}; (-40,0)*{\circlearrowright} **\dir{-};
(-37,2)*{a};
(-25,-15)*{\bullet};
(-35,-5)*{\bullet};
(-40,0)*{\circlearrowright}; (-40,10)*{\circ} **\dir{-};
(-40,10)*{\circ}; (-40,20)*{\circlearrowright} **\dir{-};
(-40,5)*{\bullet};
(-40,15)*{\bullet};
(-40,20)*{\circlearrowright}; (-50,20)*{\circ} **\dir{-};
(-50,20)*{\circ}; (-60,20)*{\circlearrowright} **\dir{-};
(-45,20)*{\bullet};
(-55,20)*{\bullet};
(-64,17)*{a^{2}};
(-60,20)*{\circlearrowright}; (-60,30)*{\circ} **\dir{-};
(-60,30)*{\circ}; (-60,40)*{\circlearrowright} **\dir{-};
(-60,25)*{\bullet};
(-60,35)*{\bullet};
(-40,20)*{\circlearrowright}; (-40,30)*{\circ} **\dir{-};
(-40,30)*{\circ}; (-40,40)*{\circlearrowright} **\dir{-};
(-40,25)*{\bullet};
(-40,35)*{\bullet};
(-36,44)*{a^{2}};
(-40,40)*{\circlearrowright}; (-50,40)*{\circ} **\dir{-};
(-50,40)*{\circ}; (-60,40)*{\circlearrowright} **\dir{-};
(-45,40)*{\bullet};
(-55,40)*{\bullet};
(-63,43)*{a^{2}};
%
%(-70,50)*{\circlearrowright}; (-65,45) **\dir{-};
%
%
%
(0,0)*{\circlearrowright}; (10,10)*{\circ} **\dir{-};
(-20,-20)*{\circlearrowright}; (-30,-30)*{\circ} **\dir{-};
(13,13)*{{\bf p}};
(-33,-33)*{{\bf p}};
(5,5)*{\bullet};
(-25,-25)*{\bullet};
(0,0)*{\circlearrowright}; (10,-10)*{\circ} **\dir{-};
(10,-10)*{\circ}; (20,-20)*{\circlearrowright} **\dir{-};
(5,-5)*{\bullet};
(15,-15)*{\bullet};
(20,-20)*{\circlearrowright}; (20,-30)*{\circ} **\dir{-};
(20,-30)*{\circ}; (20,-40)*{\circlearrowright} **\dir{-};
(20,-25)*{\bullet};
(20,-35)*{\bullet};
(20,-40)*{\circlearrowright}; (10,-50)*{\circ} **\dir{-};
(10,-50)*{\circ}; (0,-60)*{\circlearrowright} **\dir{-};
(15,-45)*{\bullet};
(5,-55)*{\bullet};
(0,0)*{\circlearrowright}; (-10,-10)*{\circ} **\dir{-};
(-10,-10)*{\circ}; (-20,-20)*{\circlearrowright} **\dir{-};
(-5,-5)*{\bullet};
(-15,-15)*{\bullet};
(-20,-20)*{\circlearrowright}; (-20,-30)*{\circ} **\dir{-};
(-20,-30)*{\circ}; (-20,-40)*{\circlearrowright} **\dir{-};
(-17,-38)*{a};
(-20,-25)*{\bullet};
(-20,-35)*{\bullet};
(-20,-40)*{\circlearrowright}; (-10,-50)*{\circ} **\dir{-};
(-10,-50)*{\circ}; (0,-60)*{\circlearrowright} **\dir{-};
(0,-56)*{a};
(-15,-45)*{\bullet};
(-5,-55)*{\bullet};
%
%(0,0)*{\circlearrowright}; (0,5)*{\bullet} **\dir{-};
%(-20,-20)*{\circlearrowright}; (-25,-15)*{\bullet} **\dir{-};
%(-20,-40)*{\circlearrowright}; (-25,-45)*{\bullet} **\dir{-};
%(0,-60)*{\circlearrowright}; (0,-65)*{\bullet} **\dir{-};
%
%
(20,-40)*{\circlearrowright}; (30,-50)*{\circ} **\dir{-};
(30,-50)*{\circ}; (40,-60)*{\circlearrowright} **\dir{-};
(40,-56)*{a};
(25,-45)*{\bullet};
(35,-55)*{\bullet};
(40,-4)*{a};
(40,0)*{\circlearrowright}; (50,-10)*{\circ} **\dir{-};
(50,-10)*{\circ}; (60,-20)*{\circlearrowright} **\dir{-};
(45,-5)*{\bullet};
(55,-15)*{\bullet};
(60,-20)*{\circlearrowright}; (60,-30)*{\circ} **\dir{-};
(60,-30)*{\circ}; (60,-40)*{\circlearrowright} **\dir{-};
(60,-25)*{\bullet};
(60,-35)*{\bullet};
(60,-40)*{\circlearrowright}; (50,-50)*{\circ} **\dir{-};
(50,-50)*{\circ}; (40,-60)*{\circlearrowright} **\dir{-};
(55,-45)*{\bullet};
(45,-55)*{\bullet};
(40,0)*{\circlearrowright}; (30,-10)*{\circ} **\dir{-};
(30,-10)*{\circ}; (20,-20)*{\circlearrowright} **\dir{-};
(35,-5)*{\bullet};
(25,-15)*{\bullet};
%
%(40,0)*{\circlearrowright}; (40,5)*{\bullet} **\dir{-};
%(40,-60)*{\circlearrowright}; (40,-65)*{\bullet} **\dir{-};
%
%(16,-22)*{a};
%(24,-22)*{a};
%(16,-38)*{a};
%(24,-38)*{a};
%(56,-22)*{a};
%(64,-22)*{a};
%(56,-38)*{a};
%(64,-38)*{a};
%(76,-22)*{a};
%(84,-22)*{a};
%(76,-38)*{a};
%(84,-38)*{a};
%
(60,-20)*{\circlearrowright}; (65,-15) **\dir{-};
(60,-40)*{\circlearrowright}; (65,-45) **\dir{-};
(70,-30)*{\ldots};
(100,0)*{\circlearrowright}; (110,-10)*{\circ} **\dir{-};
(110,-10)*{\circ}; (120,-20)*{\circlearrowright} **\dir{-};
(100,-4)*{a};
(117,-22)*{a};
(105,-5)*{\bullet};
(115,-15)*{\bullet};
(120,-20)*{\circlearrowright}; (120,-30)*{\circ} **\dir{-};
(120,-30)*{\circ}; (120,-40)*{\circlearrowright} **\dir{-};
(117,-38)*{a};
(120,-25)*{\bullet};
(120,-35)*{\bullet};
(120,-40)*{\circlearrowright}; (110,-50)*{\circ} **\dir{-};
(110,-50)*{\circ}; (100,-60)*{\circlearrowright} **\dir{-};
(100,-56)*{a};
(115,-45)*{\bullet};
(105,-55)*{\bullet};
(100,0)*{\circlearrowright}; (90,-10)*{\circ} **\dir{-};
(90,-10)*{\circ}; (80,-20)*{\circlearrowright} **\dir{-};
(95,-5)*{\bullet};
(85,-15)*{\bullet};
(80,-20)*{\circlearrowright}; (80,-30)*{\circ} **\dir{-};
(80,-30)*{\circ}; (80,-40)*{\circlearrowright} **\dir{-};
(80,-25)*{\bullet};
(80,-35)*{\bullet};
(80,-40)*{\circlearrowright}; (90,-50)*{\circ} **\dir{-};
(90,-50)*{\circ}; (100,-60)*{\circlearrowright} **\dir{-};
(85,-45)*{\bullet};
(95,-55)*{\bullet};
%
%(100,0)*{\circlearrowright}; (100,5)*{\bullet} **\dir{-};
%(120,-20)*{\circlearrowright}; (125,-15)*{\bullet} **\dir{-};
%(120,-40)*{\circlearrowright}; (125,-45)*{\bullet} **\dir{-};
%(100,-60)*{\circlearrowright}; (100,-65)*{\bullet} **\dir{-};
%
(80,-20)*{\circlearrowright}; (75,-15) **\dir{-};
(80,-40)*{\circlearrowright}; (75,-45) **\dir{-};
(0,-30)*{1};
(40,-30)*{2};
(100,-30)*{n};
(-4,0)*{a};
(4,0)*{a};
(0,4)*{a};
(0,-4)*{a};
(23,-38)*{a};
(17,-38)*{a};
(20,-44)*{a^{2}};
(63,-38)*{a};
(57,-38)*{a};
(60,-44)*{a^{2}};
(83,-38)*{a};
(77,-38)*{a};
(80,-44)*{a^{2}};
(-20,-16)*{a};
(-24,-20)*{a};
(-22,-24)*{a};
(-17,-22)*{a};
(20,-15)*{a^{2}};
(23,-22)*{a};
(17,-22)*{a};
(60,-15)*{a^{2}};
(63,-22)*{a};
(57,-22)*{a};
(80,-15)*{a^{2}};
(83,-22)*{a};
(77,-22)*{a};
(0,-64)*{a^{-1}};
(40,-64)*{a^{-1}};
(100,-64)*{a^{-1}};
(40,4)*{a^{-1}};
(100,4)*{a^{-1}};
(125,-16)*{a^{-1}};
(124,-43)*{a^{-1}};
(-24,-44)*{a^{-1}};
(-44,-4)*{a^{-1}};
(-16,24)*{a^{-1}};
(-43,23)*{a};
(-43,17)*{a};
(-37,23)*{a};
(-37,17)*{a};
(-56,24)*{a^{2}};
(-56,36)*{a^{2}};
(-44,36)*{a^{2}};
\endxy
\vspace{0.1 cm}
\begin{center}
Picture 3.
\end{center}
\vspace{0.1 cm}
\normalsize

    \par
    (Here the notations are the same as on pictures 1, 2; secondary vertices, marked with the same letter (${\bf p}$), coincide.)

    Now we show that $H_{1}'H_{2}' = G'.$ Indeed, $$(ba)^{6} \in H_{1}', \: (ba)^{5} \in H_{2}'
    \Rightarrow (ba)^{6}, \: (ba)^{5} \in H_{1}'H_{2}' \Rightarrow ba \in H_{1}'H_{2}'.$$
    Furthermore, $$(ba)^{2}(ba^{2})^{5}(ba)^{2} \in H_{2}' \Rightarrow (ba)^{2}(ba^{2})^{5}(ba)^{2}, \: ba \in H_{1}'H_{2}'
    \Rightarrow (ba^{2})^{5} \in H_{1}'H_{2}'.$$
    Moreover, $$(ba^{2})^{2} \in H_{2}' \Rightarrow (ba^{2})^{5}, \: (ba^{2})^{2} \in H_{1}'H_{2}'
    \Rightarrow ba^{2} \in H_{1}'H_{2}'.$$ Therefore, $a, b \in H_{1}'H_{2}'$, so $H_{1}'H_{2}' = G'.$ \par

    Choosing the inverse images of free generators of subgroups $H_{1}'$ and $H_{2}'$ in the same way, as in the previous case, we obtain that unimprovability of (\ref{eq3}) is proved in this case also.

    \vspace {0.2 cm}
    \par \textbf{Case 3.}  Suppose that\par
    $$ G' = \mathbb{Z}_{2}* (\mathbb{Z}_{2}\times \mathbb{Z}_{2}) \cong \langle a, b, c \:|\: a^{2}=b^{2}
    = c^{2} = bcbc = 1 \rangle $$
    \par

    We take the following subgroup as $H_{1}'\lhd G'$:
    $$H_{1}' = \langle \langle acac, \quad (ab)^{n+1} \rangle \rangle, \quad \text{then} \quad G'/H_{1}'\cong D_{n+1}\times \mathbb{Z}_{2},  $$
    where $D_{n}$ is dihedral group of order $2n$. The graph $\Psi(H_{1}')$ is represented on picture 4. \par

    \footnotesize
    \vspace{0.5 cm}
\hspace{1.6 cm}
\xy
0;/r.25 pc/:
(0,0)*{\circ}; (-5,5)*{\bullet} **\dir{-};
(-5,5)*{\bullet}; (-10,5)*{\circ} **\dir{-};
(-10,5)*{\circ}; (-15,5)*{\bullet} **\dir{-};
(-15,5)*{\bullet}; (-20,0)*{\circ} **\dir{-};
(-20,0)*{\circ}; (-25,5)*{\bullet} **\dir{-};
(0,0)*{\circ}; (5,5)*{\bullet} **\dir{-};
(5,5)*{\bullet}; (10,5)*{\circ} **\dir{-};
(10,5)*{\circ}; (15,5)*{\bullet} **\dir{-};
(15,5)*{\bullet}; (20,0)*{\circ} **\dir{-};
(20,0)*{\circ}; (25,5)*{\bullet} **\dir{-};
(0,0)*{\circ}; (-5,-5)*{\bullet} **\dir{-};
(-5,-5)*{\bullet}; (-10,-5)*{\circ} **\dir{-};
(-10,-5)*{\circ}; (-15,-5)*{\bullet} **\dir{-};
(-15,-5)*{\bullet}; (-20,0)*{\circ} **\dir{-};
(-20,0)*{\circ}; (-25,-5)*{\bullet} **\dir{-};
(0,0)*{\circ}; (5,-5)*{\bullet} **\dir{-};
(5,-5)*{\bullet}; (10,-5)*{\circ} **\dir{-};
(10,-5)*{\circ}; (15,-5)*{\bullet} **\dir{-};
(15,-5)*{\bullet}; (20,0)*{\circ} **\dir{-};
(20,0)*{\circ}; (25,-5)*{\bullet} **\dir{-};
(40,0)*{\circ}; (45,5)*{\bullet} **\dir{-};
(45,5)*{\bullet}; (50,5)*{\circ} **\dir{-};
(50,5)*{\circ}; (55,5)*{\bullet} **\dir{-};
%(55,5)*{\bullet}; (60,0)*{\circ} **\dir{-};
(40,0)*{\circ}; (35,5)*{\bullet} **\dir{-};
(40,0)*{\circ}; (45,-5)*{\bullet} **\dir{-};
(45,-5)*{\bullet}; (50,-5)*{\circ} **\dir{-};
(50,-5)*{\circ}; (55,-5)*{\bullet} **\dir{-};
%(55,-5)*{\bullet}; (60,0)*{\circ} **\dir{-};
(40,0)*{\circ}; (35,-5)*{\bullet} **\dir{-};
(30,0)*{\ldots};
(-28,5)*{{\bf p}};
(-28,-5)*{{\bf q}};
(58,5)*{{\bf p}};
(58,-5)*{{\bf q}};
(-10,7)*{a};
(-10,-7)*{a};
(10,7)*{a};
(10,-7)*{a};
(50,7)*{a};
(50,-7)*{a};
(-20,2)*{b};
(-20,-2)*{b};
(-18,0)*{c};
(-22,0)*{c};
(0,2)*{b};
(0,-2)*{b};
(2,0)*{c};
(-2,0)*{c};
(20,2)*{b};
(20,-2)*{b};
(18,0)*{c};
(22,0)*{c};
(40,2)*{b};
(40,-2)*{b};
(42,0)*{c};
(38,0)*{c};
(-10,0)*{1};
(10,0)*{2};
(25,0)*{3};
(50,0)*{n+1};
(35,0)*{n};
(-25,7)*{\ast};
\endxy
\vspace{0.2 cm}
\begin{center}
Picture 4.
\end{center}
\vspace{0.2 cm}
\normalsize
\par
(Here primary vertices are represented as $\bullet$, and secondary vertices --- as $\circ$; other notations are analogous to the notations of previous pictures. Secondary vertices, marked with the same letter (${\bf p}$ or ${\bf q}$), coincide.)
\par
    It is easy to see that $H_{1}'$ is factor-free and has finite index in $G'$. \par
     Consider the following elements of $H_{1}'$:
    \begin{equation}\label{eq20}
    w_{1}'=(acac)^{b}, \: w_{2}'=(acac)^{bab}, ..., \: w_{n}'=(acac)^{(ba)^{n-1}b}
    \end{equation}

         Notice that by choosing the maximal subtree of the graph $\Psi (H_{1}')$ in the appropriate way we can make all $n$ elements (\ref{eq20}) belong to the free generators of $H_{1}'$.

     \par
     Consider subgroup $H_{2}'\subseteq G'$, (freely) generated by the following elements:
      \begin{equation}\label{eq21}
       w_{1}', ..., w_{n}' \quad \text{from} \quad (\ref{eq20}), \quad (ab)^{acac}, \: acababa, \: abcac, \: (ac)^{(ba)^{n}b}
      \end {equation}

    The graph $\Psi(H_{2}')$ is represented on picture 5. \par

    \footnotesize
%\scriptsize
\vspace{0.5 cm}
\hspace{1.2 cm}
\xy
0;/r.25 pc/:
(-25,5)*{\bullet}; (-30,5)*{\circ} **\dir{-};
(-30,5)*{\circ}; (-35,5)*{\bullet} **\dir{-};
(-40,0)*{\circ}; (-35,5)*{\bullet} **\dir{-};
(-40,0)*{\circ}; (-45,5)*{\bullet} **\dir{-};
(-40,0)*{\circ}; (-35,-5)*{\bullet} **\dir{-};
(-40,0)*{\circ}; (-45,-5)*{\bullet} **\dir{-};
(-45,5)*{\bullet}; (-50,5)*{\circ} **\dir{-};
(-50,5)*{\circ}; (-55,5)*{\bullet} **\dir{-};
(-60,0)*{\circ}; (-55,5)*{\bullet} **\dir{-};
(-60,0)*{\circ}; (-65,5)*{\bullet} **\dir{-};
(-60,0)*{\circ}; (-55,-5)*{\bullet} **\dir{-};
(-60,0)*{\circ}; (-65,-5)*{\bullet} **\dir{-};
(-70,0)*{\circ}; (-65,5)*{\bullet} **\dir{-};
(-70,0)*{\circ}; (-65,-5)*{\bullet} **\dir{-};
(-35,-10)*{\circ}; (-25,-5)*{\bullet} **\dir{-};
(-35,-10)*{\circ}; (-45,-5)*{\bullet} **\dir{-};
(-45,-10)*{\circ}; (-35,-5)*{\bullet} **\dir{-};
(-45,-10)*{\circ}; (-55,-5)*{\bullet} **\dir{-};
(-30,7)*{a};
(-50,7)*{a};
(-72,0)*{a};
(-45,-12)*{a};
(-35,-12)*{a};
(-40,2)*{c};
(-40,-2)*{c};
(-38,0)*{b};
(-42,0)*{b};
(-60,2)*{c};
(-60,-2)*{c};
(-58,0)*{b};
(-62,0)*{b};
(0,0)*{\circ}; (-5,5)*{\bullet} **\dir{-};
(-5,5)*{\bullet}; (-10,5)*{\circ} **\dir{-};
(-10,5)*{\circ}; (-15,5)*{\bullet} **\dir{-};
(-15,5)*{\bullet}; (-20,0)*{\circ} **\dir{-};
(-20,0)*{\circ}; (-25,5)*{\bullet} **\dir{-};
(0,0)*{\circ}; (5,5)*{\bullet} **\dir{-};
(5,5)*{\bullet}; (10,5)*{\circ} **\dir{-};
(10,5)*{\circ}; (15,5)*{\bullet} **\dir{-};
(15,5)*{\bullet}; (20,0)*{\circ} **\dir{-};
(20,0)*{\circ}; (25,5)*{\bullet} **\dir{-};
(0,0)*{\circ}; (-5,-5)*{\bullet} **\dir{-};
(-5,-5)*{\bullet}; (-10,-5)*{\circ} **\dir{-};
(-10,-5)*{\circ}; (-15,-5)*{\bullet} **\dir{-};
(-15,-5)*{\bullet}; (-20,0)*{\circ} **\dir{-};
(-20,0)*{\circ}; (-25,-5)*{\bullet} **\dir{-};
(0,0)*{\circ}; (5,-5)*{\bullet} **\dir{-};
(5,-5)*{\bullet}; (10,-5)*{\circ} **\dir{-};
(10,-5)*{\circ}; (15,-5)*{\bullet} **\dir{-};
(15,-5)*{\bullet}; (20,0)*{\circ} **\dir{-};
(20,0)*{\circ}; (25,-5)*{\bullet} **\dir{-};
(40,0)*{\circ}; (45,5)*{\bullet} **\dir{-};
(45,5)*{\bullet}; (50,0)*{\circ} **\dir{-};
%(50,5)*{\circ}; (55,5)*{\bullet} **\dir{-};
%(55,5)*{\bullet}; (60,0)*{\circ} **\dir{-};
(40,0)*{\circ}; (35,5)*{\bullet} **\dir{-};
(40,0)*{\circ}; (45,-5)*{\bullet} **\dir{-};
(45,-5)*{\bullet}; (50,0)*{\circ} **\dir{-};
%(50,-5)*{\circ}; (55,-5)*{\bullet} **\dir{-};
%(55,-5)*{\bullet}; (60,0)*{\circ} **\dir{-};
(40,0)*{\circ}; (35,-5)*{\bullet} **\dir{-};
(30,0)*{\ldots};
%(-28,5)*{p};
%(-28,-5)*{q};
%(58,5)*{p};
%(58,-5)*{q};
%
(-10,7)*{a};
(-10,-7)*{a};
(10,7)*{a};
(10,-7)*{a};
(52,0)*{a};
(-20,2)*{b};
(-20,-2)*{b};
(-18,0)*{c};
(-22,0)*{c};
(0,2)*{b};
(0,-2)*{b};
(2,0)*{c};
(-2,0)*{c};
(20,2)*{b};
(20,-2)*{b};
(18,0)*{c};
(22,0)*{c};
(40,2)*{b};
(40,-2)*{b};
(42,0)*{c};
(38,0)*{c};
(-10,0)*{1};
(10,0)*{2};
(25,0)*{3};
(35,0)*{n};
%(50,0)*{n+1};
%
(-25,7)*{\ast};
\endxy
\vspace{0.2 cm}
\begin{center}
Picture 5.
\end{center}
\vspace{0.2 cm}
\normalsize

    It is easy to see that $H_{2}'$ is factor-free and has finite index in $G'$. \par

    Now we show that $H_{1}'H_{2}' = G'$. Indeed,
    $$acac \in H_{1}', \: (ab)^{acac} \in H_{2}' \Rightarrow acac, \: (ab)^{acac} \in H_{1}'H_{2}' \Rightarrow ab \in H_{1}'H_{2}'. $$
    Furthermore, $$abcac, \: ab \in H_{1}'H_{2}' \Rightarrow cac, \: acac \in H_{1}'H_{2}' \Rightarrow a, b \in H_{1}'H_{2}'. $$
    Moreover, $$acababa, a, b \in H_{1}'H_{2}' \Rightarrow c \in H_{1}'H_{2}', $$
    so $H_{1}'H_{2}' = G'$.
    \par
    Choosing the inverse images of free generators of subgroups $H_{1}'$ and $H_{2}'$ in the same way, as in the first case (only here we take $w_{1}', ..., w_{n}'$ from (\ref{eq20})), we obtain that unimprovability of (\ref{eq3}) is proved in this case also.

    \vspace {0.2 cm}
    \par \textbf{Case 4.}  Consider the last case. Suppose that \par
    $$ G' = \mathbb{Z}_{2} * \mathbb{Z}_{2} * \mathbb{Z}_{2} \cong \langle a, b, c \:|\: a^{2}=b^{2}
    =c^{2} = 1 \rangle $$
    \par
     Let $G_{*}', H_{1*}', H_{2*}'$ denote groups $G', H_{1}', H_{2}'$ from the previous case respectively. Then $$G_{*}' = \mathbb{Z}_{2} * (\mathbb{Z}_{2} \times \mathbb{Z}_{2}) \cong G'/  \langle \langle bcbc \rangle \rangle $$
    Consider a subgroup $$ H_{1}' \lhd G', H_{1}' = \langle \langle bcbc, acac, (ab)^{n+1} \rangle \rangle,$$
    As the generators of $H_{2}'$ we take the union of the elements from (\ref{eq21}) and the following elements:
    \begin{equation}\label{eq22}
    (bcbc)^{r_{j}}, \quad r_{j} =  1, ba, (ba)^{2}, ..., (ba)^{n}, a, aca.
    \end{equation}
    Then $$G'/H_{1}' \cong D_{n+1}\times \mathbb{Z}_{2} \cong G_{*}'/H_{1*}',$$
    and it is easy to see the correspondence between $\Psi(H_{1*}'), \Psi(H_{2*}')$ (graphs of subgroups of $G_{*}'$) and $\Psi(H_{1}'), \Psi(H_{2}')$ (graphs of subgroups of $G'$) respectively. \par
    It is also easy to see that subgroups $H_{1}'$ and $H_{2}'$ are factor-free and have finite index in $G'$. The same arguments as in the previous case prove that $H_{1}'H_{2}' = G'$. \par
    Choosing the inverse images of free generators of subgroups $H_{1}'$ and $H_{2}'$ in the same way, as in the previous cases, we obtain that unimprovability of (\ref{eq3}) is proved in this case also, and that means theorem 2 is proved.

    \vspace {0.3 cm}
    The author thanks A.A.Klyachko for many useful remarks.

    \vspace{0.3 cm}
    Faculty of Mechanics and Mathematics, Moscow State University, Russia \par
    {\it E-mail address:} zakhar.sasha@gmail.com \par \par
    This work was supported by the Russian Foundation for Basic Research, project no. 11-01-00945.

\end{document}